\documentclass{amsart}

\usepackage[numbers]{natbib}
\usepackage{amsmath, amssymb,stmaryrd}
\usepackage{type1cm}
\usepackage{tensor}
\usepackage{enumerate}
\usepackage{graphicx}
\usepackage[czech,english]{babel}
\usepackage{mdwtab}
\usepackage[OT1]{fontenc}
\usepackage[all]{xy}

\newcommand{\bbbn}{\mathbb{N}}

\usepackage{color}
\definecolor{shadecolor}{gray}{.85}%
\definecolor{tintedcolor}{gray}{.80}%
\usepackage{framed}%
\definecolor{mytintedcolor}{gray}{.95}%
\makeatletter
\newdimen\svparindent
\newcounter{tmpthm}
\newenvironment{mytinted}{%
  \MakeFramed {\FrameRestore}}%
{\endMakeFramed}
{\endlist\end{mytinted}\egroup}
\makeatother
\newcommand{\cat}[1]{\mbox{\,\bf #1}}
\newcounter{mydots}
\def\newdot{\addtocounter{mydots}{1}\noindent{\bf \themydots.~}}

\newtheorem{theorem}{Theorem}

\newtheorem{lemma}[theorem]{Lemma}

\theoremstyle{definition}

\theoremstyle{remark}

\begin{document}

\title{Towards a Characterization of Universal Categories}
\author{Jaroslav Ne{\v s}et{\v r}il}
\address{Jaroslav Ne{\v s}et{\v r}il\\
Computer Science Institute of Charles University (IUUK and ITI)\\
   Malostransk\' e n\' am.25, 11800 Praha 1, Czech Republic}
\email{nesetril@iuuk.mff.cuni.cz}
\thanks{Supported by grant ERCCZ LL-1201 
and CE-ITI P202/12/G061, and by the European Associated Laboratory ``Structures in
Combinatorics'' (LEA STRUCO)}

\author{Patrice Ossona~de~Mendez}
\address{Patrice~Ossona~de~Mendez\\
Centre d'Analyse et de Math\'ematiques Sociales (CNRS, UMR 8557)\\
  190-198 avenue de France, 75013 Paris, France
	--- and ---
Computer Science Institute of Charles University (IUUK)\\
   Malostransk\' e n\' am.25, 11800 Praha 1, Czech Republic}
\email{pom@ehess.fr}
\thanks{Supported by grant ERCCZ LL-1201 and by the European Associated Laboratory ``Structures in
Combinatorics'' (LEA STRUCO), and partially supported by ANR project Stint under reference ANR-13-BS02-0007}

\date{\today}
\begin{abstract}
In this note we characterize, within the framework of the theory of finite set, those categories of graphs that are {\em algebraic universal} in the sense that every concrete category embeds in them.
The proof of the characterization is based on the sparse--dense dichotomy and its model theoretic equivalent.
\end{abstract}
\maketitle
\section{Introduction}
\label{sec:intro}
A category $\mathcal K$ is {\em algebraic universal} if every concrete category embeds in it. The name comes  from examples:
algebraic universal categories include simple algebraic structures as well as the class of all graphs (sets with one binary relation).
Algebraic universal categories have been the subject of intensive studies
\cite{Adamek2004,PT}. Particularly, many subcategories of the category of graphs were shown to be algebraic universal, too \cite{adams1981,babai1978high,Hell1973,PT}. 

The aim of this note is to provide a characterization of those subcategories of the category of graphs that are universal. Unexpectedly this is related to (and in fact coincides with) the characterization of somewhere dense classes of graphs. All these notions will be introduced in Section~\ref{sec:prel}. 

At this place let us remark that we deal only with finite graphs and categories induced by them, so this paper is in fact written in the theory of finite sets (so $\bbbn$ is a proper class here), see Section~\ref{sec:comments}.

The main result of this note is the following:
Denote by $\cat{Gra}$ the category of all finite undirected graphs, and by
$\overrightarrow{\cat{Gra}}$ the category of all finite oriented graphs. Recall that an {\em oriented graph} is a directed graph in which at most one arc exists between any two vertices.

\begin{theorem}
\label{thm:main}
	For a monotone subcategory $\mathcal K$ of $\cat{Gra}$ the following three statements are equivalent:
	\begin{enumerate}
		\item There exists a subcategory $\overline{\mathcal K}$ of $\overrightarrow{\cat{Gra}}$, each member of which is an orientation of a member of $\mathcal K$,  which embeds the category $\overrightarrow{\cat{Gra}}$.
		\item There exists a subcategory $\overline{\mathcal K}$ of $\overrightarrow{\cat{Gra}}$, each member of which is an orientation of a member of $\mathcal K$,  which embeds the simplicial category $\Delta$.
		\item $\mathcal K$ is somewhere dense.
	\end{enumerate}
\end{theorem}

In other words, some class of orientations of graphs from $\mathcal K$ is universal (in the theory of finite sets) if and only if the class $\mathcal K$ is somewhere dense.

This in turn leads to a new, high level, algebraic characterization of somewhere dense classes. Yet another one in the already long list, see \cite{Sparsity, SurveyND}.

The paper is organised as follows: In Section~\ref{sec:prel} we recall all the relevant notions and put the universality question in the context of category theory related to concreteness and representation of posets.
In Section~\ref{sec:proof} we prove the main result by a combination of model theory and combinatorial methods. In Section~\ref{sec:comments} we recast the problem of universality in the context of the classification of sparse classes of graphs, and display a perhaps surprising gap in the descriptive complexity of classes representing groups, monoids, and categories.

\section{Preliminaries}
\label{sec:prel}
First we recall several notions of category theory.

A category is {\em concrete} if it is isomorphic to a subcategory of the category $\cat{Set}$ of sets and mappings. A necessary condition for a category to be concrete is {\em Isbell condition}  \cite{Isbell1963}.
This condition was proved to be sufficient by Freyd \cite{freyd1973} and,
through the explicit construction of a faithful functor to $\cat{Set}$, by Vin{\'a}rek \cite{vinarek1976} .   Vin{\'a}rek's construction
has, moreover, the following property: for countable categories with finite
sets of morphisms between fixed objects, the functor has finite values. Thus  Freyds' theorem holds also in the finite set theory.  Precisely, if the considered class  is countable and the set of homomorphisms between any two objects is finite,  then the class is isomorphic to a subcategory of the category  of finite sets if and only if Isbell's condition holds.
See \cite{HN} for a concise description of these results.

In this context, another interesting result is Ku{\v{c}}era's theorem \cite{kuvcera1971}, which asserts that every category is a factorization of a concrete one (like classes of homotopy equivalent maps, which was the original motivation of \cite{freyd1973}). 
Also this theorem holds  in its finite set theory version.
 
 As a culmination of researches by Prague category group in the sixties, it has been proved that the category \cat{Gra} of all graphs (finite or infinite) with homomorphisms between them is  algebraic universal for all concrete categories. Explicitly, for every concrete category $\mathcal K$ there is an embedding of $\mathcal K$ into \cat{Gra}.

These results led to an intensive research, and various subcategories of \cat{Gra} were shown to be algebraic universal.

The basic techniques used in these proofs was model-theoretical first-order interpretation, then called {\v s}{\'\i}p, indicator, or replacement construction.

It is perhaps surprising that in this paper we can provide a characterization of monotone subcategories of \cat{Gra} that are algebraic universal. Here {\em monotone} means that the class (of graphs) is closed under taking (non necessarily induced) subgraphs. 

As this paper deals with finite models we restrict from now in the setting of  finite set theory, thus to finite graphs and to  embedding into the category of finite graphs.

In order to formulate our main results we have to recall the basics of the nowhere dense--somewhere dense dichotomy. For a comprehensive treatment, see e.g. \cite{Sparsity} or \cite{SurveyND}.

Somewhere dense classes of graphs were introduced by the authors in \cite{ND_logic,ND_characterization}. Recall that a class of graphs $\mathcal C$ is {\em nowhere dense} if, for every integer $p$ 
there exists an integer $N(p)$ such that the $p$-th subdivision 
${\rm Sub}_p(K_{N(p)})$ of $K_{N(p)}$ is a subgraph of no graph in $\mathcal C$, and the class $\mathcal{C}$ is {\em somewhere dense}, otherwise. So,  a monotone class $\mathcal C$ is nowhere dense if an only if there exists $N:\bbbn\rightarrow\bbbn$ such that for every $p\in\bbbn$ it holds
${\rm Sub}_p(K_{N(p)})\notin\mathcal{C}$. Nowhere dense classes found various applications in designing fast (almost linear) algorithms \cite{POMNII,Grohe2013,Dvovr'ak2011}. Particular cases of nowhere dense classes are {\em classes with bounded expansion} \cite{POMNI}. These are characterized by the property that for every integer $p$ 
there exists an integer $N(p)$ such that no $p$-th subdivision of a graph with minimum degree at least $N(p)$ is a subgraph of a graph in $\mathcal C$. Such classes have strong structural and algorithmic properties \cite{Sparsity}.

Let $\mathcal C$ be a class of structures of a fixed signature. A first-order formula $\phi(\mathbf x,\mathbf y)$ is said to have the {\em order property} with respect to $\mathcal C$ if it has the {\em $n$-order property} for all $n$, i.e. if for every $n$ there exist a structure $\mathbf M\in\mathcal C$ and tuples $\mathbf a_0,\dots,\mathbf a_{n-1},b_0,\dots,\mathbf b_{n-1}$ of elements of $\mathbf M$ such that $\mathbf M\models \phi(\mathbf a_i,\mathbf b_j)$ holds if and only if $i<j$. A class $\mathcal C$ of structures is called {\em stable} if there is no such formula with respect to $\mathcal C$. It is easy to see that $\mathcal C$ is stable if and only if there is no formula $\psi(\mathbf u,\mathbf v)$  with $|\mathbf u|=|\mathbf v|$, such that for every $n$ there exist a structure $M\in \mathcal C$ and tuples $\mathbf c_0,\dots,\mathbf c_{n-1}$ of elements of $\mathbf M$ such that $\mathbf M\models \psi(\mathbf c_i,\mathbf c_j)$  holds if and only if $i< j$, i.e. $\psi$ orders the tuples linearly.
Stability and the ($n$-)order property come from stability theory \cite{Shelah1990,Ensley1996}, where they are defined for the class of models of a complete first-order theory. 
In \cite{Adler2013}, Adler and Adler prove the following theorem (see also \cite{SurveyND}).
\begin{theorem}
\label{thm:Adler}
	Let $\mathcal C$ be a monotone class of coloured digraphs of a fixed finite signature, and let $\underline{\mathcal C}$ be the class of the underlying undirected graphs. The following conditions are equivalent.
	\begin{enumerate}
		\item  $\underline{\mathcal C}$ is nowhere dense;
		\item ${\mathcal C}$ is stable;
	\end{enumerate}
\end{theorem}

This interplay of model theoretic and combinatorial notions is the key to our main result.
\section{Characterization}
\label{sec:proof}
In view of the context of our main result (outlined in Section~\ref{sec:intro}) it suffices to prove the following two results.

\begin{lemma}
\label{lem:1}
	 Let $\mathcal D$ be a monotone somewhere dense class of  undirected graphs. Then there exists a class $\mathcal C$ of oriented graphs, each member of which is an orientation of a graph in $\mathcal D$, which represents the  category $\overrightarrow{\cat{Gra}}$ of oriented graphs.
\end{lemma}
	\begin{proof}
			Let $d$ be such that $\mathcal C$ contains the $d$-subdivision of every complete graph $K_n$ (here we use the assumption that $\mathcal C$ is monotone).  Let $(I,a,b)$ be the circuit of length $3(d+1)$, where vertices $a,b$ are  linked by a directed path (from $a$ to $b$) of length $d+1$.
			For a given oriented  graph $\vec G$, denote by $\vec G\ast (I,a,b)$ the directed graph which arises from $\vec G$ by replacing every arc $(u,v)$ of $\vec G$ by a copy of $(I,a,b)$ in such a way that $a$ is identified to $u$ and $b$ to $v$ (all other vertices in distinct copies being distinct).

The only circuits of  $G\ast (I,a,b)$ with length at most  $3(d+1)$ occur as copies of $(I,a,b)$.  It follows that any homomorphism 
$f:\vec G\ast(I,a,b)\rightarrow \vec H\ast(I,ab)$ is induced in a unique way by a homomorphism $g:\vec G\rightarrow \vec H$.
\end{proof}

\begin{lemma}
\label{lem:2}
If a class of directed graphs represents the 	the simplicial category $\Delta$, then it is somewhere dense.
\end{lemma}
\begin{proof}
Let $\mathcal C$ be a class of directed graphs that represents the 	the category $\Delta$.
	For sake of simplicity, we assume that every directed graph $\vec G$  we consider has vertex set $0,1,2,\dots,|\vec{G}|-1$, and we denote by $E(x,y)$ the relation expressing the existence of an arc from $x$ to $y$.

Then there is a functor $\Phi$, mapping each ordinal $[n]=\{0,1,\dots,n\}$
 to  a directed graph $\Phi([n])\in\mathcal C$, and 
bijectively mapping order preserving maps $f:[i]\rightarrow [j]$ into homomorphisms $\Phi(f):\Phi([i])\rightarrow\Phi([j])$ in such a way that $\Phi(f\circ g)=\Phi(f)\circ\Phi(g)$. 

Let $\vec G_n=\Phi([n])$, let $a=|\vec G_0|-1$ and $b=|\vec{G}_1|-1$. We define the formula
	$$
	\nu(x_0,\dots,x_a):=\bigwedge_{\vec{G}_0\models E(i,j)} E(x_i,x_j),
	$$
	which asserts that $i\mapsto x_i$ is a homomorphism from $\vec{G}_0$.

	There are exactly two order preserving maps from $[0]$ to $[1]$, namely  $f_s:0\mapsto 0$, and $f_t:0\mapsto 1$.
	Let $\phi_s=\Phi(f_s)$ and $\phi_t=\Phi(f_t)$.
	
	Then we define
	\begin{multline*}
		\eta(x_0,\dots,x_a,y_0,\dots,y_a):=\\
		(\exists z_0\dots z_b)\quad
		\Bigl[
		\bigwedge_{i=0}^a (x_i=z_{\phi_s(i)})
		\wedge 		\bigwedge_{i=0}^a (y_i=z_{\phi_t(i)})
		\wedge        \bigwedge_{\vec{G}_1\models E(i,j)} E(z_i,z_j)
		\Bigr]
	\end{multline*}

The meaning of formula $\eta$ is as follows: $\vec{G}\models\eta(x_0,\dots,x_a,y_0,\dots,y_a)$  expresses  that there exist 
a homomorphism $h:\vec{G}_1\rightarrow\vec{G}$ and homomorphisms 
$g_s,g_t:\vec{G}_0\rightarrow\vec{G}$, such that 
$g_s(i)=x_i$, $g_t(i)=y_i$, $g_s=h\circ \phi_s$, and $g_t=h\circ \phi_t$.
In other words, naming $z_i=h(i)$, there exist $z_0,\dots,z_b$ such that
$i\mapsto z_i$ is a homomorphism $\vec G_1\rightarrow\vec G$,  $x_i=z_{\phi_s(i)}$,  and $y_i=z_{\phi_t(i)}$. (Note that $\phi_s$ and $\phi_t$ are known to be homomorphisms.)

Let $n$ be an ordinal. There are exactly $n+1$ order preserving maps $g_i:[0]\rightarrow[n]$, that are naturally ordered in such a way that for every $i,j\in[n]$ it holds $i<j$ if and only if $g_i(0)<g_j(0)$. 
 In other words,
 for every two order preserving maps $g_i,g_j:[0]\rightarrow[n]$ there exists an order preserving map $h:[1]\rightarrow[n]$ such that $g_i=h\circ f_s$ and $g_j=h\circ f_t$ if and only if $i<j$. 
It follows that 
  for every two homomorphisms $\hat g,\hat g':\vec{G}_0\rightarrow\vec{G}_n$ there exists an homomorphism $\hat h:\vec{G}_1\rightarrow\vec{G}_n$ such that $\hat g=\hat h\circ \phi_s$ and $\hat g=\hat h\circ \phi_t$ if and only if the (uniquely determined) integers $i,j$ such that $\hat g=\Phi(g_i)$ and $\hat g'=\Phi(g_j)$ are such that  $i<j$. 
  
  Define the $n+1$ tuples $\overline{x}^i=(x_0^i,\dots,x_a^i)$ by 
  $x_i^j=\Phi(g_j)(i)$. In other words, let $i\mapsto x_i^j$ be the homomorphism $\Phi(g_j)$. Then the above properties rewrites as
  $$\vec{G}_n\models\eta(\overline{x}^i,\overline{x}^j)\quad\iff\quad i<j.$$
  It follows that $\mathcal C$ has the order property hence, by Theorem~\ref{thm:Adler}, is somewhere dense.
\end{proof}

\section{Comments}
\label{sec:comments}
\newdot
 Let us add	few remarks putting the results of this work in a broader context.
Representation of categories were first investigated in the special cases of groups, monoids, and small categories. This line of research directly relates to our main result.

For groups, the representation can be done by graphs \cite{frucht1939}, and even by $3$-regular graphs \cite{sabidussi1957graphs}. However this cannot be done by geometrically restricted graphs, like planar graphs or, more generally, by any proper minor closed class of graphs \cite{Babai1974}.

For monoids, the representation can be done by graphs with arbitrary girth (this is also possible by the above construction) but not by $3$-regular or even $k$-regular graphs (for any fixed $k$). In fact Babai and Pultr \cite{BabPul} 
showed that any class of graphs which represents all finite monoids has to contain a subdivision of any complete graph. However, using large girth representations and using characterization of classes with bounded expansion \cite{POMNI, Sparsity}, one can easily see that finite monoids can be represented by graphs in a bounded expansion class $\mathcal C_0$. In particular, one can put $\mathcal C_0$ to be the class of all graphs of the form $G\ast(\vec C_{2n},a,b)$, where $\vec C_{2n}$ is a circuit  of length $2n$, where $n$ is the order of $G$.

Consider small categories (in the theory of finite sets, that is finite categories). Let us enumerate all non-isomorphic small categories of graphs as $\mathcal K^1,\mathcal K^2,\dots,\mathcal K^n,\dots$.

Let $\mathcal K^i$ have objects $G_1^i,\dots,G_{t(i)}^i$. The category 
 $\mathcal K^i$  will be represented by oriented graphs of the form
 $G_j^i\ast(\vec C_{2N_i},a,b)$, where $N_i\geq \sum_{j=1}^{t(i)}|V(G_j^i)|$ and $N_i<N_j$ whenever $i<j$. On sees easily that the class $\mathcal C_1$ of all such graphs $G_j^i\ast(\vec C_{2N_i},a,b)$ has bounded expansion: for any fixed integer $d$ and any graph $H$, if the $d$-th subdivision of $H$ is a subgraph of a graph in $\mathcal C_1$ then $H$ is $2$-degenerate with possibly finitely many exceptions (derived from graphs $G_j^i\ast(\vec C_{2N_i},a,b)$ for small $i$). The class $\mathcal C_1$ represents all the small categories $\mathcal K^1,\mathcal K^2,\dots,\mathcal K^n,\dots$ by an application of Cayley-MacLane representation.

However to represent arbitrary categories (in the theory of finite sets) we have to jump over nowhere dense classes, right to somewhere dense classes. This descriptive complexity gap is surprising.

It would be interesting to find a more direct combinatorial proof of the fact that representing special categories leads to bounded subdivisions of arbitrarily large complete graphs. Such examples of groups and monoids we found in 
\cite{Babai1974,BabPul}.


\newdot
In this context one should note that the representation of posets and thin categories can be achieved by oriented paths, trees, or outerplanar graphs \cite{Hubivcka2005,HubN}.

Let us summarize these facts in a schematic table.
\vspace{2mm}

\newcommand{\mc}{\multicolumn{1}}
\begin{tabular}{|l||p{.3\textwidth}|p{.3\textwidth}|}
\mc{l}{}	&\mc{c}{$+$}&\mc{c}{$-$}\\
\hlx{v[0,1,2,3]hv}
Posets&	oriented trees, cycles, or paths&undirected bipartite \\
\hlx{vhv}
Groups& bounded degree&proper minor closed\\
\hlx{vhv}
Monoids& bounded expansion & proper topological minor closed\\
\hlx{vhv}
Small categories& bounded expansion&proper topological minor closed\\
\hlx{vhv}
Concrete categories&somewhere dense&nowhere dense\\
\hlx{vh}
\end{tabular}
\vspace{2mm}

\newdot We restricted ourselves to the finite set theory (i.e. to finite graphs). The situation for infinite graphs and categories is less clear. 
On the other hand most examples of special algebraic universal categories are obtained from some basic examples (like the category of graphs) 
by first-order  interpretation (like replacement operation in the above proof of Lemma~\ref{lem:2}). As the basic examples contain complete graphs of any size this leads then to $p$-subdivisions of large complete graphs. The main result of this paper shows that this is necessarily so. 

\providecommand{\bysame}{\leavevmode\hbox to3em{\hrulefill}\thinspace}
\providecommand{\MR}{\relax\ifhmode\unskip\space\fi MR }
\providecommand{\MRhref}[2]{%
  \href{http://www.ams.org/mathscinet-getitem?mr=#1}{#2}
}
\providecommand{\href}[2]{#2}

\end{document}